\documentclass[12pt]{article}
\newcommand{\qed}{\hfill\rule{4pt}{8pt}\par\vspace{\baselineskip}}
\input amssymb.sty
\input mac2.sty
\def\T{{\mathbb{T}}}
\newtheorem{thm}{Theorem}[section]
\newtheorem{defn}[thm]{Definition}

\newtheorem{exs}[thm]{Examples}
\newtheorem{lem}[thm]{Lemma}
\newtheorem{pro}[thm]{Proposition}
\newtheorem{cor}[thm]{Corollary}
\newtheorem{rem}[thm]{Remark}

\def\bpol{{^{\displaystyle \diamond}}}
\def\pol{{^{\displaystyle \circ}}}

\begin{document}
\title{On the continuity of separately continuous bihomomorphisms}
\author {R. Beattie and H.-P. Butzmann}
\date{}
\maketitle
\begin{abstract}
Separately continuous bihomomorphisms on a product of convergence or
topological groups occur with great frequency. Of course, in
general, these need not be jointly continuous. In this paper,
we exhibit some results of Banach-Steinhaus type and use these to
derive joint continuity from separate continuity. The setting of
convergence groups offers two advantages. First, the continuous
convergence structure is a powerful tool in many duality arguments.
Second, local compactness and first countability, the usual
requirements for joint continuity, are available in much greater
abundance for convergence groups.\sp

\end{abstract}

\section{Introduction}

Let G, H and L be topological groups and $u: G \times H \to L$ a
separately continuous bihomomorphism. One need look no further than
the evaluation mapping $\o : \G_{co} \times G \to \T$ to see that
such bihomomorphisms need not be jointly continuous \cite{MP}. The
problem of determining conditions on G, H and L so that u is jointly
continuous is a difficult one and has a long history in the
literature. In fact the problem has been studied extensively in the
larger context of topological spaces. Over several decades, many
results have appeared (see e.g. \cite{N}, \cite{P1},
\cite{P2},\cite{P3} or \cite{HT}) guaranteeing points of joint
(quasi-)continuity for various combinations of topological spaces G,
H and L. Recurring themes were the notions of compactness and
countability, the latter usually appearing in some form of the Baire
property. The topological group case simplifies the problem
considerably since joint (quasi-)continuity at one point implies
joint continuity.

In this paper we address this problem in the context of convergence
groups. Apart from providing greater generality than topological
groups, this permits the use of continuous convergence in duality
arguments. With the aid of the notion of a $g$-barrelled group, we
establish theorems of Banach-Steinhaus type and use these together
with duality arguments to establish the joint continuity of
bihomomorphisms. The main result is the following: If G is a
$g$-barrelled group, H a locally compact convergence group and L a
locally quasi-convex topological group, every separately continuous
bihomomorphism $u: G \times H \to L$ is jointly continuous. The
generality of this result can be seen as various special cases
recover many of the results in the literature.

Let $X$ be a set and suppose that to each $x$ in $X$ is associated a
collection $\l(x)$ of filters on $X$ satisfying for all $x \in X$:

\noindent (i) the ultrafilter $\dot{x}:=\{ A \subseteq X: x \in A \}
\in \l (x)$,\sp (ii) if $\cF \in \l(x)$ and $\cG \in \l(x)$, then
$\cF \cap \cG \in \l(x)$,\sp (iii) if $\cF \in \l(x)$, then $\cG \in
\l(x)$ for all filters $\cG \supseteq \cF$.

The totality $\lambda$ of filters $\lambda (x)$ for $x$ in $X$ is
called a {\bf convergence structure} for $X$, the pair $(X,
\lambda)$ a {\bf convergence space} and filters $\cal F$ in $\lambda
(x)$ {\bf convergent} to $x$. A convergence space $(X, \lambda)$
will usually be denoted by $X$ if no confusion arises. We write
${\cal F} \rightarrow x$ instead of ${\cal F} \in \lambda (x).$  A
mapping $f: X \rightarrow Y$ between the two convergence spaces $X$
and $Y$ is {\bf continuous} if $f({\cal F}) \rightarrow f(x)$ in $Y$
whenever ${\cal F} \rightarrow x$ in $X$.

Let $G$ be a group (all groups will be assumed to be Abelian)
and assume $\lambda$ is a convergence structure on $G$. The pair
$(G, \lambda)$ is a {\bf convergence group} if $\lambda$ is
compatible with the group operations, i.e., if the mapping
$$-: G \times G \to G,\quad (x,y) \mapsto x-y$$
is continuous. This means that if ${\cal F} \rightarrow x$ and
${\cal G} \rightarrow y$ in $G$, then the filter ${\cal F} -  {\cal
G}$ generated by $\{ A-B \mid A \in {\cal F}, B \in {\cal G} \}$
converges to $x - y$ in $G$.

Every topological space is a convergence space, the convergent
filters at any point being precisely those finer than the
neighbourhood filter.  Likewise, every topological group is a
convergence group. The converse statements fail. Convergence groups
need not be topological.

A convergence space $X$ is called {\bf Hausdorff} if limits are
unique, i.e., if ${\cal F} \rightarrow p$ and ${\cal F} \rightarrow
q$ in $X$, then $p =q.$ It is called {\bf compact} if each
ultrafilter converges and {\bf locally compact} if it is Hausdorff and
each convergent filter contains a compact set.

Let $G, H$ be convergence groups and $\Gamma(G, H)$ the space of
continuous group homomorphisms from $G$ to $H$. The {\bf continuous
convergence} structure on $\Gamma (G, H)$ is the coarsest
convergence structure on $\Gamma (G, H)$ making the evaluation
mapping
$$\omega: \Gamma(G, H) \times G \to H, \st (\f,x) \mapsto \f(x)$$
continuous. A filter $\Phi \rightarrow \f$ in $\Gamma_{c}(G, H)$ if,
whenever $\cF \rightarrow x$ in $G$, the filter $\o(\Phi \times
\cF)$ converges to $\o(\f,x) = \f(x)$ in $H$. The continuous
convergence structure is compatible with the group $\Gamma (G, H)$
and the resulting convergence group is denoted $ \Gamma_{c} (G, H)$.
When $H = \T = \R/\Z$, one obtains the continuous dual
$\Gamma_{c}G$, the canonical dual space of a convergence group $G$.
Note that, when $G$ is a topological group, the continuous dual
$\Gamma_{c} (G)$ is a locally compact convergence group. In general
it is not topological, but this is so if $G$ is  locally
compact. In this case the continuous convergence structure is the
compact-open topology.

If $u : G \to H$ is a continuous homomorphism between convergence
groups, then $u^* : \G H \to \G G$ is defined by $u^*(\psi) = \psi
\circ u$. It is continuous if both character groups are either endowed
with the continuous convergence structure or the weak topology
(defined below). In this way $\G_c$ becomes a functor which has strong categorical
properties. It is a left adjoint and takes final structures to
initial structures, in particular quotients to embeddings and direct
limits to  inverse limits.

If $G$ is any convergence group then the canonical mapping $\k_G : G
\rightarrow \G_c\G_cG$ defined by
$$\k_G(x)(\f) = \f(x) \st \mbox{ for all } x \in G \mbox{ and all } \f \in \G G$$
is always continuous. A convergence group $ G$ is called {\bf
embedded} if $\k_{G}$ is an isomorphism onto its range and {\bf
reflexive} if $\k_G$ is an isomorphism.

The {\bf weak topology} on $\G(G,H)$ is the initial topology induced
by the family of mappings $(\f \mapsto \f(x))_{x \in G}$. The
resulting topological group is denoted by $\G_s(G,H)$. As above, when
$H = \T$, this becomes  $\G_sG$.

Finally, $\rho : \R \to \T$ denotes the canonical projection and we
set
$$\T_+ = \rho([-1/4, 1/4])$$
If $\T$ is realized as the unit circle, this is the right half of it.

Further information on convergence spaces and in particular
convergence groups can be found in \cite{Bi} and \cite{BB}.

\section{$g$-barrelled convergence groups}

In a linear setting, topological vector spaces and convergence
vector spaces, the Banach-Steinhaus Theorem relates pointwise
bounded and equicontinuous sets as well as pointwise and
continuously convergent sequences (see e.g. \cite{Bou}, \cite{BB},
\cite{BS}). Whereas the notion of equicontinuity generalizes very
naturally to the setting of convergence groups, the notion of
(pointwise) boundedness is usually not available and must be
replaced.
\begin{defn}

Let $G,H$ be convergence groups. A set $M \subseteq \Gamma (G,H)$ is
called {\bf equicontinuous} if and only if, for all filters $\cF$
which converge to $0$ in $G$, the filter $M(\cF)$ converges to $0$
in $H$. Here $M(\cF)$ denotes the filter generated by $\{M(F) : F
\in M \} = \{\o(M \times F) : F \in \cF\}$.
 \end{defn}

It is clear that, when $G$ and $H$ are topological groups, this
coincides with the usual definition of equicontinuity.

As the next proposition shows, equicontinuity is preserved as $G$
and $H$ pass to final and initial structures respectively.

\begin{pro}\label{equi wrt initial and final} Let $G$ and $H$ be
convergence groups. If $G$ carries the final group convergence
structure with respect to a family of homomorphisms $(u_i : G_i \to
G)_{i \in I}$ and $H$ carries the inital group convergence structure
with respect to family of homomorphisms $(v_j : H \to H_j)_{i \in
I}$ then a set $M \subseteq \G(G,H)$ is equicontinuous if and only
if for all $i \in I$ and $j \in J$ the set
$$v_j \circ M \circ u_i = \{v_j \circ w \circ u_i : w \in M\}$$
is an equicontinuous subset of $\G(G_i, H_j)$.
\end{pro}

\pf An easy argument shows that $v_j \circ M \circ u_i$ is
equicontinuous for all  $i \in I$ and $j \in J$ if $M$ is
equicontinuous. To show the converse, assume that $\cF \to 0 \in G$. Since $G$ carries the final group
convergence structure with respect to $(u_i)$ there are $i_1,
\ldots, i_n \in I$ and filters $\cF_k \to 0 \in G_{i_k}$ such that
$$\cF \supseteq u_{i_1}(\cF_1) + \cdots + u_{i_n}(\cF_n)$$
By assumption, $v_j \circ M \circ u_{i_k}(\cF_k)$ converges to $0$
in $H_j$ for all $j \in J$ and $k \in \{1, \ldots, n\}$ and
therefore $v_j(M(\cF)) = v_j \circ M(\cF)$ converges to $0$ for all
$j$. Since $H$ carries the initial group convergence structure with
respect to $(v_j)$ we get $M(\cF) \to 0 \in H$ as desired. \qed\\

What makes equicontinuous sets valuable for our purposes is the
following result (see \cite[2.4.2]{BB} for a general formulation).

\begin{pro}\label{s=c} Let $G$ and $H$ be convergence groups and let $M \subseteq
\G(G,H)$ be an equicontinuous set. Then the weak topology and the
continuous convergence structure coincide on $M$.
\end{pro}

The following notion was defined for topological groups by E.
Martin-Peinador and V. Tarieladze in \cite{MPT} and \cite{CMPT}.

\begin{defn} A convergence group $G$ is called {\bf $g$-barrelled}
if the compact subsets of $\Gamma_s(G)$ are equicontinuous.
\end{defn}

The standard examples of $g$-barrelled topological groups are
countably $\check{C}$ech-complete topological groups, so, in
particular, complete metrizable or locally compact ones. Also
separable Baire or metrizable hereditarily Baire groups are
$g$-barrelled (\cite{N}, \cite{T}, \cite{CMPT}). Finally, the
additive group of a barrelled topological vector space is
$g$-barrelled (\cite{MPT}). To obtain an example of a non-topological $g$-barrelled
convergence group we recall that a topological group $G$ is said to
{\bf respect compactness} if each
 $\s(G,\G G)$-compact subset of $G$ is compact.

\begin{pro}\label{g-barrelled dual} If $G$ is a reflexive topological
group that respects compactness, then $\G_cG$ is $g$-barrelled.
\end{pro}

\pf Since $G$ is reflexive, the natural mapping $\k_G : G \to
\G_c\G_cG$ is an isomorphism and therefore $\k_G : (G, \s(G,\G G))
\to \G_s\G_cG$ is an isomorphism. If $M \subseteq \G_s\G_cG$ is
compact, then $\k_G^{-1}(M)$ is a $\s(G, \G G)$-compact subset of
$G$ and therefore compact. This implies that $M =
\k_G(\k_G^{-1}(M))$ is compact. So it is equicontinuous by
the Arzel\`a-Ascoli-Theorem (\cite[2.5.6]{BB}).\qed

\begin{cor}If $G$ is a nuclear group, then $\G_cG$ is $g$-barrelled.
\end{cor}

\pf If $G$ is a complete nuclear group, it is reflexive by
\cite[8.4.19]{BB} and it respects compactness by \cite{BMP}. Therefore
$\G_cG$ is $g$-barrelled by \ref{g-barrelled dual}. If $G$ is an
arbitrary nuclear group then its
completion $\widetilde{G}$ is nuclear by \cite[21.4]{Au}. Also
$\G_c(G) = \G_c(\widetilde{G})$ by \cite[8.4.4]{BB} and so the result
follows. \qed

It should be noted that the reflexive locally convex topological
vector spaces which respect compactness
are precisely the Montel spaces \cite[Theorem 1.4]{RTA}.

The next several propositions derive permanence properties of
$g$-barrelled convergence groups.

\begin{pro}\label{final}  \hspace*{0em}\sp
(i) Let $G$ and $G'$ be convergence groups with the same underlying
group such that $\G G = \G G'$. If the identity mapping $id : G \to
G'$ is continuous, then $G$ is $g$-barrelled if $G'$ is.\sp
(ii) A convergence  group which carries the final group convergence
structure with respect to a family of group homomorphisms from
$g$-barrelled convergence groups is $g$-barrelled.\sp
(iii) A topological group which carries the final group topology with
respect to a family of group homomorphisms from $g$-barrelled
topological groups is $g$-barrelled.
\end{pro}

\pf
 (i) Evidently $\G_sG = \G_sG'$, so if  $M \subseteq \G_s G$ is
compact, then $M$ is compact in $\G_sG'$ and therefore
equicontinuous. So if $\cF$ converges to $0$ in $G$ then it
converges to $0$ in $G'$ and therefore $M(\cF)$ converges to $0$. So
 $M$ is an equicontinuous subset of $\G G$.\sp
(ii) Assume that  $G$ carries the final group convergence structure
with respect to a family of group homomorphisms $(G_i \to G)_{i \in
I}$ such that all $G_i$ are $g$-barrelled. If $\cF$ converges to $0$
in $G$ there are finitely many $i_1, \ldots, i_n \in I$ and filters
$\cF_j$ converging to zero in $G_{i_j}$ such that
$$\cF \supseteq u_{i_1}(\cF_1) + \cdots + u_{i_n}(\cF_n)$$
Take any compact subset $M$ of  $\G_sG$. Then $u_i^*(M)$ is compact
in $\G_s(G_i)$ for all $i$ and therefore equicontinuous. So
$M(u_{i_j}(\cF_j)) = u_{i_j}^*(M)(\cF_j)$ converges to $0$ in $\T$
and so does
$$M(u_{i_1}(\cF_1)) + \cdots + M(u_{i_n}(\cF_n))$$
The claim now follows from
$$M(\cF) \supseteq M(u_{i_1}(\cF_1) + \cdots + u_{i_n}(\cF_n))
\supseteq M(u_{i_1}(\cF_1)) + \cdots + M(u_{i_n}(\cF_n))$$
(iii) This is  \cite[1.9]{CMPT}.\qed


\begin{lem} \label{prep products} Let $(G_i)_{i \in I}$ be a family of
  convergence groups,
$G = \prod_{i \in I}G_i$ their product and let $e_i : G_i \to G$ be the
natural injections. If $M \subseteq \G_s(G)$ is compact, then there
is a finite subset $I_0 \subseteq I$ such that $\f \circ e_i = 0$
for all $\f \in M$ and all $i \in I \setminus I_0$.
\end{lem}

\pf Since the mapping
$$\G_cG \lra \bigoplus_{i \in I} \G G, \st \f \mapsto (\f \circ e_i)$$
is an isomorphism by \cite[8.1.18]{BB}, we will regard the elements of $\G G$ as elements
in $\bigoplus \G G_i$. The claim then is that there is a finite set
$I_0 \subseteq I$ such that $\f_i = 0$ for all $\f \in M$ and all $i
\in I \setminus I_0$. So assume that this is not true. For shorter
reference, for all $\f \in \G G$ we set $C(\f) = \{i \in I : \f_i
\not=0\}$. These sets are all finite. Define inductively sequences
$(\f_n)$ in $M$ and $(i_n)$ in $I$ in the following way:

Choose any $\f_1 \in M, \; \f_1 \not=0$ and any $i_1 \in I$ such
that $\f_{1,i_1} \not=0$. ( Here and in what follows $\f_{n,i}$
denotes the $i$-th component of $\f_n$.) Assume that $\f_1, \ldots,
\f_{n-1}$ and $i_1, \ldots i_{n-1}$ have been chosen. Then there is
a $\f_n \in M$ such that $C(\f_n) \not\subseteq C(\f_1) \cup \ldots
\cup C(\f_{n-1})$. Choose any $i_n \in C(\f_n) \setminus (C(\f_1)
\cup \ldots \cup
C(\f_{n-1}))$. Then $\f_{n,i_n} \not=0$.\\
Note that, by construction, we have
$$\f_{j,i_n} = 0 \st \mbox{ for all } j < n$$
and therefore, in particular, $i_j \not=i_n$ for all $j < n$.\\
Now choose any $x = (x_i) \in G$ such that $x_i = 0$ if $i \notin
\{i_n : n \in \N\}$. Then for all $r \in \N$ we get:
$$\f_r(x) = \sum_{i \in I}\f_{r,i}(x_i) = \sum_{n \in \N}\f_{r,
i_n}(x_{i_n}) = \sum_{r \ge n} \f_{r, i_n}(x_{i_n}) = \sum_{n<r}
\f_{r, i_n}(x_{i_n}) + \f_{r,i_r}(x_{i_r})$$

Set $\T_0
= \rho([-1/16,1/16])$. We show that for each finite set $J \subseteq
I$ there is an element $x \in G$ such that $x_i = 0$ for all $i \in
J$
and $\f_n(x) \notin \T_0$ for all but finitely many $n$.\\
Choose a finite set $J \subseteq I$ and a $k \in \N$ such that $i_n
\notin J$ for all $n \ge k$. Now define $x \in G$ in the following
way: $x_i = 0$ if $i \notin \{i_n : n \in \N\}$ and also $x_i = 0$
for all $i \in \{i_n : n < k\}$. Then $x_i = 0$ if $i \in J$. Define
$x_{i_n}$ for all $n \ge k$ inductively as follows: One has
$$\f_k(x) = \f_{k, i_k}(x_{i_k})$$
and since $\f_{k,i_k} \not=0$ there is some $x_{i_k} \in G_{i_k}$
such that $\f_{k, i_k}(x_{i_k}) \notin \T_0$.\\
If $x_{i_k}, \ldots, x_{i_{r-1}}$ have been constructed, we get
$$\f_r(x) = \sum_{n<r} \f_{r, i_n}(x_{i_n}) + \f_{r,i_r}(x_{i_r})$$
If $\sum_{n<r} \f_{r, i_n}(x_{i_n}) \notin\T_0$, then set $x_{i_r} =
0$ otherwise there is some $x_{i,r} \in G_{i_r}$ such that
$\f_{r,i_r}(x_{i_r}) \notin \T_+$ and then $\f_r(x) \notin \T_0$.

Assume now that the sequence $(\f_n)$ has a cluster point $\psi \in
\G(G)$. Then there is a finite set $J \subseteq I$ such that $\psi_i
= 0$ for all $i \in I \setminus J$. Choose $x$ as above. Then
$\psi(x) = 0$ and so there must be infinitely many $n$ such that
$\f_n(x) = \f_n(x) - \psi(x) \in \T_0$, contradicting the
construction of $x$.\qed

\begin{pro} Let $(G_i)_{i \in I}$ be a family of $g$-barrelled
convergence groups. Then $\prod_{i \in I}G_i$ is $g$-barrelled.
\end{pro}

\pf Set $G = \prod_{i \in I}G_i$ and let $M \subseteq \G_s(G)$ be a
compact set. Since $e_i^* : \G_s(G) \to \G_s(G_i)$ is continuous for
all $i$, also $e_i^*(M)$ is compact in $\G_s(G_i)$ and therefore
 equicontinuous.
 By Lemma \ref{prep products}, there are elements $i_1, \ldots, i_n \in I$ such that
 $ \f \circ e_i = 0$ for all $\f \in M$ and all $i \not= i_1, \ldots, i_n$. Take any filter
 $\cF$ which converges to $0$ in $G$, then $p_i(\cF)$ converges to
 $0$ in $G_i$, where $p_i$ denotes the projection, and therefore
 $M(e_{i}(\pi_{i}(\cF))) = e_{i}^*(M)(\pi_{i}(\cF))$ converges
 to $0$ for all $i$. Choose a zero
 neighbourhood $U$ in $\T$, then there is  a zero neighbourhood $V$ in $\T$ such
 that $nV = V + \cdots + V \subseteq U$. Then there is a set $F \in \cF$
 such that
 $$M(e_{i_1}(\pi_{i_1}(F))) + \cdots + M(e_{i_n}(\pi_{i_n}(F)))
 \subseteq nV  \subseteq U$$
Take now any $\f \in M$ and $x  \in F$. Then we have
$$\f(x) = \sum_{i \in I}e_i^*(\f)(x_i) =  \sum_{i \in I}\f \circ e_i(x_i)
= \sum_{j =1}^n\f \circ e_{i_j}(x_{i_j}) = \sum_{j
=1}^n\f(e_{i_j}(\pi_{i_j}(x))) \in U$$ and so $M(F) \subseteq U$.
\qed

Locally quasi-convex topological groups will be of  importance in the
  sequel and so we introduce them here as well as the locally
  quasi-convex modification.

A subset $A$ of a topological group $G$ is called {\bf
quasi-convex} if for each $x \in G \setminus A$ there is a
character $\f \in \G G$ such that $\f(A) \subseteq \T_+$ while
$\f(x) \notin \T_+$.  Furthermore,
$G$ is called {\bf locally quasi-convex} if it has a zero
neigbourhood base consisting of quasi-convex sets. As it turns out
each Hausdorff topological group $G$ is locally quasi-convex and
Hausdorff if and only if it is embedded (see \cite[8.4.7]{BB}.

If $G$ is a convergence group, then the finest locally
quasi-convex topology on $G$ which is coarser than the
convergence structure of $G$ is called the {\bf locally
quasi-convex modification} of $G$ and the resulting topological group
is denoted by $\tau(G)$.  In order to give an explicit
description thereof,  for
subsets $A \subseteq G$ and  $H\subseteq \G G$ we define
$$A\pol = \{\f \in \G G : \f(A) \subseteq \T_+\}$$
and
$$H\bpol = \{x \in G : H(x) \subseteq \T_+\}$$

In this terminology $A$ is quasi-convex if and only if $A = A\pol \bpol $.

\begin{thm} Let $G$ be a convergence group. Then
$$\cB :=\{H\bpol : H \subseteq \G G \mbox{ equicontinuous}\}$$
is a zero neighbourhood base of the locally quasi-convex modification
of $G$.
\end{thm}

\pf Clearly $\cB$ is a filter basis consisting of symmetric sets
and, if $H$ is an equicontinuous subset
of $\G G$ containing $0$, then $H + H$ is also equicontinuous and
$$(H + H)\bpol + (H +H) \bpol \subseteq H\bpol$$
Therefore $\cB$ is the zero neighbourhood basis of a locally
quasi-convex topology $\tau$ on $G$. If $\cF$ converges to $0$ in $G$, then
$H(\cF)$ converges to $0$ in $\T$ and so there is some $F \in \cF$
such that $H(F) \subseteq \T_+$. This gives $F \subseteq H\bpol$ and
so the zero neighbourhood filter of $\tau$ is contained in $\cF$ which gives
the continuity of the identity mapping $id : G \to (G,
\tau)$. Finally, if $\mu$ is any locally quasi-convex topology on $G$
coarser then that of $G$ and $V$ is any quasi-convex zero neighbourhood
in $(G, \mu)$ then $V\pol$ is an equicontinuous subset of $\G (G,
\mu)$ and therefore of $\G G$. Consequently, $V = V\pol\bpol \in \cB$ and so
$id : G \to (G, \mu)$ is continuous. \qed

\begin{pro}If $G$ is a convergence group then $\G G$ and $\G \tau(G)$
  share the same equicontinuous subsets.
\end{pro}

\pf Clearly each equicontinuous subset if $\G \tau(G)$ is equicontinuous
in $\G G$. One the other hand, if $H$ is an equicontinuous subset of
$\G G$ then $H\bpol$ is a zero neighbourhood of $\tau(G)$ and so
$H\bpol\pol$ is an equicontinuous subset of $\G \tau(G)$ containing
$H$. \qed

\begin{cor} \label{G and tauG} A convergence group $G$ is $g$-barrelled if and only if
  $\tau(G)$ is.
\end{cor}

\begin{cor} A topological group which carries the final locally
  quasi-convex group topology with
respect to a family of group homomorphisms from $g$-barrelled
topological groups is $g$-barrelled.
\end{cor}

\pf This follows from \ref{final}(ii) and \ref{G and tauG}. \qed\\

The concept of $g$-barrelledness allows us to relate the compact
subsets of $\G_s(G,H)$ and the equicontinuous subsets of
$\G(G,H)$. The following two theorems can be thought of as theorems
of Banach-Steinhaus type.

\begin{thm} \label{target=dual} Let $G$ and $H$ be convergence groups.
If G is $g$-barrelled and  $H$ is locally compact then each compact
subset of $\Gamma_s (G,\Gamma_c (H))$ is equicontinuous.
\end{thm}

{\bf Proof.} Let $M$ be a compact subset of $\Gamma_s (G,
\Gamma_c(H))$ and assume that ${\cal F} \to 0$ in $G$. We have to
show that $M({\cal F}) \to 0$ in $\Gamma_c(H)$. So let ${\cal H} \to
z$ in $H$. Since $H$ is locally compact, $\cal H$ contains a compact
set $K$ and $M({\cal F})({\cal H})$ is finer than $M({\cal F})(K).$
We claim that $M({\cal F})(K)$ converges to $0$ in $\T$ which will
give the desired result.\\
 Consider the mapping
 $$ T : \Gamma_s (G,\Gamma_c(H)) \times H \to \Gamma_s(G)$$
 given by $T(u,y)(x) = u(x)(y)$ for all
 $(u,y) \in  \Gamma(G,\Gamma_c(H)) \times H$ and all  $x \in G$.
 An easy calculation shows that $T$ is continuous and so
$T(M \times K)$ is compact in $\Gamma_s(G)$. Since $G$ is
$g$-barrelled, $T(M \times K)$ is equicontinuous in $\Gamma(G)$.
Hence $M({\cal F})(K) = T(M \times K)(\cF) \to 0$ in $\T$ as
required.\qed

\begin{thm}\label{compact implies equi} Let $G$ be a $g$-barrelled
convergence group and $L$ a Hausdorff locally quasi-convex
topological group. Then the compact subsets of   $\G_s(G,L)$ are
equicontinuous.
\end{thm}

{\bf Proof.} Since $L$ is a Hausdorff locally quasi-convex
topological group, it is an embedded convergence group, and
therefore isomorphic to a subgroup of $\Gamma_c \Gamma_c L$. Set $H
= \G_cL$. Then $H$ is locally compact. So if $M$ is a compact subset
of $\G_s(G,L)$, it can be considered a compact subset of $\G_s(G,
\G_cH)$ and is therefore equicontinuous by Proposition
\ref{target=dual}. Clearly $M$ is then equicontinuous in $\G(G,L)$.
\qed\\

Since the weak topology and the continuous convergence structure
coincide on equicontinuous sets by Proposition \ref{s=c},
Theorem \ref{compact implies equi} gives conditions under which each
compact subset of $\G_s(G,L)$ is even compact in $\G_c(G,L)$.

\section{Joint Continuity of bihomomorphisms}

In this section we make use of the results of the previous section
to extract the joint continuity of separately continuous
bihomomorphisms in several special cases. A key observation here is
the following:

\begin{pro} \label{translation} Let $G, H$ and $L$ be convergence groups and
$u : G \times H \to L$ be a separately continuous bihomomorphism.
Then the mapping
$$u_H : H \lra  \G_s(G,L)$$
defined by $u_H(y)(x) = u(x,y)$ is continuous. Furthermore, $u$ is
jointly continuous if and only if
$$u_H : H \lra  \G_c(G,L)$$
is continuous.
\end{pro}

\pf The first part is clear. Now from the universal property of the
continuous convergence, $u_H$ is continuous if and only if the
mapping
$$\o \circ (id_G \times u_H) : G \times H \to L$$
is continuous. But evidently $\o \circ (id_G \times u_H) = u$ and so
the proof follows. \qed

\begin{pro} \label{joint under hyp}  Let $G, H$ and $L$ be convergence groups
 such that the compact subsets of $\Gamma_s (G,L)$ are equicontinuous. Assume
 further that $u : G \times H \to L$ is a separately continuous bihomomorphism.
 Then $u$ is jointly continuous in either of the following two cases:\sp
(i) $H$ is locally compact.\sp
(ii) $G$ and $H$ are first countable and $L$ is topological.
\end{pro}

{\bf Proof.}  By
\ref{translation} we must show that $u_H: H \to \Gamma_c(G,L)$ is
continuous. \sp
(i)  $u_H : H \to \G_s(G,L)$ is continuous by \ref{translation}. Let $\cF \to y_0$ in $H$. Then $u_H(\cF)$
converges to $u_H(y_0)$ in $\G_s(G,L)$. Since $\cF$ contains a
compact set $u_H(\cF)$ contains a compact subset of $\G_s(G,L)$. By
assumption, this set is equicontinuous and so $u_H(\cF)$ converges
to $u_H(y_0)$ in $\G_c(G,L)$ by \ref{s=c}. \sp (ii) We first show
that $u_H$ is sequentially continuous: If $(y_n)$ is a sequence
which converges to $y_0$ in $H$, then $B = \{y_n : n \in \N\} \cup
\{y_0\}$ is compact subset of $H$ and therefore $u(B)$ is compact
and hence equicontinuous. Again, this implies that $(u_H(y_n))$
converges to
$u_H(y_0)$ in $\G_c(G,L)$, and so $u_H$ is sequentially continuous.\\
If now $(x_n,y_n)$ is a sequence in $G \times H$ which converges to
$(x_0,y_0)$ in $G \times H$ then $(u_H(y_n))$ converges to
$u_H(y_0)$ in $\G_c(G,L)$ by the first part and so $(u(x_n,y_n)) =
(u_H(y_n)(x_n))$ converges to $u(x_0, y_0) = u_H(y_0)(x_0)$ in $H$.
Since $G$ and $H$ are first countable the claim follows. \qed\\

 From
\ref{translation}, \ref{joint under hyp} and \ref{compact implies
equi} we get the main result of this section:

\begin{thm} \label{main4} Let $G$ and $H$ be convergence groups, $G$
$g$-barrelled, and let $L$ be a locally quasi-convex topological
group.
 Then every separately
continuous bihomomorphism $u: G \times H \to L$ is jointly
continuous in either of the following cases:\sp
(i) $H$ is locally compact.\sp
(ii) $G$ and $H$ are first countable.
\end{thm}

 Part (ii) of the above theorem yields joint continuity
 results for first countable convergence groups. It uses duality
 arguments. One can also obtain joint continuity in first countable
 situations using standard Baire category techniques.

\begin{pro} Let $G$ and $L$ be topological groups, $G$ Baire, and let $H$
a first countable convergence group. Then each separately continuous
bihomomorphism $u : G \times H \to L$ is jointly continuous.
\end{pro}

\pf Since $u$ is separately continuous, it suffices to show that $u$
is continuous at $(0,0)$. So assume that $\cF$ converges to $0$ in
$H$. Since $H$ is first countable there is a filter $\cV \subseteq
\cF$ with a countable base $(V_n)$ which also converges to $0$. Take
a closed zero neighbourhood $W$ in $L$.  For all $n \in \N$ consider
the set
$$A_n = \{x \in G : u(x \times V_n) \subseteq W\}$$
We first claim that $\bigcup A_n = G$: Take any $x \in G$, then
$u(x,0) = 0$. Since $u(x, \cdot)$ is continuous, there is an $n \in
\N$ such that $u(x \times V_n) \subseteq W$. This means $x \in A_n$.
Next, we show that each $A_n$ is closed. For all $y \in V_n$ we have
$u(A_n \times y) \subseteq W$. Since $u$ is separately continuous,
we have
$$u(\ol{A_n} \times y) \subseteq \ol{u(A_n \times y)} \subseteq \ol{W} = W$$
and so $u(x,y) \in W$ for all $x \in \ol{A_n}$ and all $y \in V_n$.
This gives $u(x \times V_n) \subseteq W$ for all $x \in \ol{A_n}$
and therefore $\ol{A_n}\subseteq A_n$.\sp Since $G$ is a Baire space,
some $A_k$ has an interior point $x_0$ and so there is a zero
neighbourhood $U$ in $G$ such that $x_0+ U \subseteq A_k$. This
gives
$$u(x_0 + x,y) \in W \st \mbox{ for all } x \in U, \, y \in V_k$$
Furthermore, there is a zero neighbourhood $V$ in $H$ such that
$u(x_0 \times V) \subseteq W$ and so
$$u(x_0, y) \in W \st \mbox{ for all } y \in V$$
Finally we get, for all $x \in U$ and $y \in V_k \cap V$,
$$u(x,y) = u(x_0 + x,y) - u(x_0,y) \in W - W$$
which shows that $u$ is continuous at $(0,0)$. \qed\\

 One factor which makes the results of Theorem \ref{main4} strong is
the size of the class of $g$-barrelled convergence groups. Even if one restricts
oneself to topological groups, this class remains large. As seen in
the previous section, it includes all countably
$\check{C}$ech-complete topological groups and is closed under the formation of
arbitrary products and inductive limits.

\begin{exs} \hspace*{0em}\sp
(i) Let $G$ be an  inductive limit of locally compact topological groups, $H$ a
locally compact convergence group
and  $L$  a locally quasi-convex topological group. Then every separately
continuous bihomomorphism $u: G \times H \to L$ is jointly
continuous.\sp
(ii) Let $G$  be a convergence  inductive
limit of complete metrizable topological groups, $H$ a metrizable
topological group  and $L$ a locally
quasi-convex topological group. Then every
separately continuous bihomomorphism $u: G \times H \to L$ is
jointly continuous.
\end{exs}

It should be mentioned that, in general, inductive limits depend very heavily on the setting
in which they are taken. It is a consequence of Proposition
\ref{final} and \ref{G and tauG}, however, that the convergence group
(or topological group or locally quasi-convex group) inductive limit of $g$-barrelled groups is also
$g$-barrelled. Thus in Examples 3.5(i), $G$ may be any appropriate inductive limit.

Another factor which adds scope to the results of Theorems 3.3 and 3.4 is the fact that, for convergence groups,
the notions of local compactness and first countability are not as restrictive as for topological groups.
Convergence group inductive limits preserve both properties. Also, for any topological group $G$,
the continuous character group $\Gamma_cG$ is locally compact.

\begin{exs}\hspace*{0em}\sp
(i) Let $G$ be  a convergence inductive limit of
 Baire groups, $H$ a convergence inductive limit of metrizable
 topological groups and
  $L$  any topological group. Then every separately
continuous bihomomorphism $u: G \times H \to L$ is jointly
continuous.\sp
(ii) Let $G$ be a convergence inductive limit of complete metrizable
 topological groups,
 $H$ and $L$ topological groups, $L$ locally quasi-convex.
Then every separately continuous bihomomorphism $u: G \times
\Gamma_c (H )\to L$ is jointly continuous.\sp
(iii) Let $G$ and $H$ be separable metrizable topological groups, $G$ complete  and $L$ any topological group.
Then any separately continuous bihomomorphism $u: G
\times \Gamma_c (H) \to L$ is jointly continuous.\sp
(iv) Let $G, H$ and $L$ be topological groups, $G$ nuclear and $L$
 locally quasi-convex.  Then every separately continuous bihomomorphism $u: \Gamma_c
(G) \times \Gamma_c (H )\to L$ is jointly continuous.
\end{exs}

\begin{rem} {\rm
Results on joint continuity can be viewed as results on triples $(G,H,L)$ of convergence
or topological groups. In such situations relaxing restrictions on one variable often
requires tightening them on another. Consider the following:\sp
(i) In \cite{HT} it is shown that separate continuity implies joint continuity  if $G$ and $H$ are both countably
$\check{C}$ech complete and $L$ is metrizable. This is a relaxation of
the condition of the local compactness of $H$ in \ref{main4}(i) but is
much more restrictive on $G$ and $L$.\sp
(ii) One can easily generalize the notion of sequential barrelledness
defined in \cite{CMPT} and \cite{MPT} to convergence groups.
A convergence group $G$ is sequentially barrelled if every convergent sequence in $\Gamma_s(G)$ is
equicontinuous. This is a large class of groups which includes all $g$-barrelled groups and all Baire
groups. It is possible to imitate the proof of Theorem \ref{main4} to obtain joint continuity for $G$ first countable
and sequentially barrelled, $H$ first countable and $L$
 a second countable locally quasi-convex topological group. This is a relaxation of the conditions
on $G$ but is much more restrictive on $L$.\sp
(iii) If $G$ is assumed only to be $g$-barrelled, it does not appear that one can relax the condition
of local compactness on $H$ very far. If $G$ is a complete metrizable
topological group and $H = \Gamma_{co}(G)$ is its
Pontryagin dual, then $G$ is $g$-barrelled and $H$ is a $k$-space and $k$-group \cite{Ch}, but the evaluation mapping
$\omega: G \times H \to \T$ is not jointly continuous.
}
\end{rem}

R. Beattie, Department of Mathematics and Computer Science, Mount
Allison University Sackville, N.B., Canada, E4L 1E6
(rbeattie@mta.ca).

H.-P. Butzmann, Fakult\"{a}t f\"{u}r Mathematik und Informatik,
Universit\"{a}t Mannheim, Mannheim, Germany,
(butzmann@math.uni-mannheim.de).

\end{document}